\documentclass{article}
\usepackage{amsthm}
\usepackage[T1]{fontenc}
\usepackage[frenchb, english]{babel}
\usepackage[hmargin=4cm,vmargin=4cm]{geometry}
\usepackage{amssymb}
\usepackage{amsmath}
\usepackage{pstricks}

\def\aro{{A}^{\!\!\!\raise5pt\hbox{$\scriptstyle \circ$}}}
\def\aroi{{A}^{\!\!\!\raise4pt\hbox{$\scriptscriptstyle \circ$}}}
\def\cro{\smash{{C}^{\!\!\!\raise5pt\hbox{$\scriptstyle \circ$}}}}
\def\croi{\smash{{C}^{\!\!\!\raise4pt\hbox{$\scriptscriptstyle \circ$}}}}

\begin{document}

\newtheorem{Th}[subsubsection]{Théorème}
\newtheorem{Pro}[subsubsection]{Proposition}
\newtheorem{De}[subsubsection]{Définition}
\newtheorem{Prt}[subsubsection]{Propriété}
\newtheorem{Prts}[subsubsection]{Propriétés}
\newtheorem{Le}[subsubsection]{Lemme}
\newtheorem{Hyp}[subsubsection]{Hypothèse}
\newtheorem{Cor}[subsubsection]{Corollaire}



\title{A linear version of Dawson-Gärtner's theorem}
\date{March 21, 2011}
\author{Pierre Petit\\ \\
Universit\'e Paris Sud}
\maketitle

\selectlanguage{english}

\begin{abstract}
We prove a linear version of Dawson-Gärtner's theorem saying that weak large deviations principles and the equivalence of ensembles $s = -p^*$ are preserved through linear projective limits. This version takes into account the linear structure underlying Cramér's theorem and turns out to be well appropriate to it.
\end{abstract}

\selectlanguage{french}

\begin{abstract}
Nous démontrons une version linéaire du théorème de Dawson-Gärtner assurant que les principes de grandes déviations faibles et l'équivalence d'ensemble $s = -p^*$ sont conservés par passage aux limites projectives d'espaces vectoriels. Cette version tient compte de la structure linéaire sous-jacente au théorème de Cramér, pour lequel elle s'avère bien adaptée.
\end{abstract}

\medskip

\textbf{Remarque :} La prochaine version de ce texte sera en langue anglaise.

\section{Introduction}

L'objectif de ce texte est de donner les définitions les plus naturelles de la pression et de l'entropie, et de faire apparaître le lien entre les deux fonctions, autrement dit les fondements de l'équivalence d'ensembles. Le but n'est pas ici de montrer l'équivalence d'ensembles, ce qui ne peut se faire que moyennant des hypothèses sur la suite $(\mu_n)_{n \geqslant 1}$ des mesures considérées : par exemple, lois des moyennes empiriques $\overline{X}_n$ d'une suite $(X_n)_{n \geqslant 1}$ de variables i.i.d. Mais nous ébauchons un cadre pour la théorie de Cramér qui peut être vu comme l'analogue du cadre de Varadhan \cite{Var66} pour les principes de grandes déviations (PGD) forts.

Dans le cadre de Varadhan, fortement inspiré du théorème de Sanov, la compacité joue un rôle particulier, notamment dans les bonnes fonctions de taux qui apparaissent. Nous affaiblissons ici le rôle de la compacité en exploitant le rôle de la linéarité dans la théorie de Cramér, et en particulier pour la question de l'équivalence d'ensembles $s = -p^*$. On aboutit à une version linéaire du \emph{théorème de Dawson-Gärtner} assurant le transport des PGD faibles et de l'égalité $s = -p^*$ par limite projective d'espaces vectoriels. Ce résultat est très adapté au cadre des champs asymptotiquement découplés développé dans \cite{CeP11} qui fait naturellement apparaître une structure de limite projective, et, dans le cas indépendant, il permet d'obtenir l'égalité $s = -p^*$ dans un cadre très général (cf. \cite{Pet11a}). On notera qu'en revanche la borne supérieure pour les convexes ne passe pas aux limites projectives, de même que l'égalité duale $p = (-s)^*$. Ce texte doit beaucoup à \cite{LeP95} pour les notations utilisées et pour plusieurs résultats montrés au passage (comme la caractérisation des PGD faibles).

Nous avons découpé l'exposition en quatre parties, ajoutant à chaque étape les hypothèses nécessaires : nous introduisons d'abord la pression, puis l'entropie en lien avec la notion de PGD faible (au passage, nous énonçons aussi un \emph{principe de contraction} linéaire) ; les deux dernières parties sont consacrées au transport, par limite projective, des PGD faibles et de l'égalité $s = -p^*$.

\section{Pression}

\subsection{Cadre et définitions}

Soient $X$ un ensemble et $\mathcal{F}$ une tribu sur $X$. Soient $(\mu_n)_{n \geqslant 1}$ une suite de mesures de probabilité sur $\mathcal{F}$ et $(v_n)_{n \geqslant 1}$ une suite de réels strictement positifs. On note $\mathcal{M}$ l'ensemble des fonctions mesurables de $X$ dans $[-\infty , +\infty]$ et on définit, pour tout $\varphi \in \mathcal{M}$,
$$
\underline{p}(\varphi) := \liminf_{n \to \infty} \frac{1}{v_n} \log \int e^{v_n \varphi} d\mu_n \qquad \textrm{et} \qquad \overline{p}(\varphi) := \limsup_{n \to \infty} \frac{1}{v_n} \log \int e^{v_n \varphi} d\mu_n
$$
On appelle \emph{pression} de la suite $(\mu_n, v_n)_{n \geqslant 1}$ la fonction $p = \overline{p}$.

\subsection{Convexité de la pression}

On étend l'addition à $[-\infty , +\infty ]$ via
$$
\forall a \in \mathbb{R} \qquad a \,\textup{\d{$+$}}\, (\pm\infty) = \pm\infty \qquad \textrm{et} \qquad (-\infty) \,\textup{\d{$+$}}\, (+\infty) = -\infty
$$
(on gardera, sauf pour $(-\infty) \,\textup{\d{$+$}}\, (+\infty)$, la notation $+$) ainsi que la multiplication par un réel via
$$
\begin{array}{ll} \forall \alpha > 0 & \alpha \cdot (\pm\infty) = \pm\infty \\
\forall \alpha < 0 & \alpha \cdot (\pm\infty) = \mp\infty
\end{array} \qquad \textrm{et} \qquad 0 \cdot (\pm\infty) = 0
$$
On munit $\mathcal{M}$ de l'addition et de la multiplication par un réel terme à terme. On dit qu'une fonction $g : \mathcal{M} \rightarrow [-\infty , +\infty]$ est \emph{convexe} si
$$
\forall (\varphi, \psi) \in X^2 \quad \forall \alpha \in [0, 1] \qquad g(\alpha \varphi + (1-\alpha)\psi) \leqslant \alpha g(\varphi) \,\textup{\d{$+$}}\, (1-\alpha) g(\psi)
$$

\begin{Pro}
La pression est une fonction convexe.
\end{Pro}

\textbf{Démonstration :} Pour tout $n \geqslant 1$, l'inégalité de Hölder permet de montrer que
$$
\varphi \in \mathcal{M} \mapsto \frac{1}{v_n} \log \int e^{v_n \varphi} d\mu_n \in [-\infty , +\infty]
$$
est convexe. On conclut en remarquant qu'une limite supérieure de fonctions convexes est convexe.\qed

\section{Entropie et principes de grandes déviations}

\subsection{Cadre et définitions}

Soient $X$ un ensemble, $\mathcal{F}$ une tribu sur $X$ et $\tau$ une topologie séparée\footnote{Toutes les topologies seront supposées séparées, on ne le précisera plus. La raison est d'ordre culturel et pratique... Notons toutefois que cette hypothèse est inutile : il suffit de remplacer chaque occurrence de \og{}compact\fg{} par \og{}quasi-compact\fg{}.} sur $X$. Soient $(\mu_n)_{n \geqslant 1}$ une suite de mesures de probabilité sur $\mathcal{F}$ et $(v_n)_{n \geqslant 1}$ une suite de réels strictement positifs. On définit, pour tout $x \in X$,
$$
\underline{s}(x) := \inf_{\substack{A \in \mathcal{F}\\ \aroi \ni x}}  \liminf_{n \to \infty} \frac{1}{v_n} \log \mu_n(A) \qquad \textrm{et} \qquad \overline{s}(x) := \inf_{\substack{A \in \mathcal{F}\\ \aroi \ni x}} \limsup_{n \to \infty} \frac{1}{v_n} \log \mu_n(A)
$$
Bien entendu, $\underline{s} \leqslant \overline{s}$. On appelle\footnote{On dit, d'habitude, l'entropie de la suite $(\mu_n)_{n \geqslant 1}$, de vitesse $(v_n)_{n \geqslant 1}$, relativement à $\tau$. C'est par souci de concision que nous l'écrivons ainsi. Notons aussi que \cite{LeP95} réserve le terme d'\emph{entropie} au cas où $\underline{s} = \overline{s}$ est concave.} \emph{entropie} de la suite $(\mu_n, v_n)_{n \geqslant 1}$, relativement à $\tau$, la fonction $s = \underline{s}$. Par construction, l'entropie est la plus grande fonction vérifiant la borne inférieure :

\medskip

\textsf{(BI)} pour tout $A \in \mathcal{F}$,
$$
\liminf_{n \to \infty} \frac{1}{v_n} \log \mu_n(A) \geqslant \sup_{\aroi} s
$$

On dit que $(\mu_n, v_n)_{n \geqslant 1}$ vérifie un \emph{principe de grandes déviations} (PGD) si la borne supérieure suivante est vérifiée :

\medskip

\textsf{(BS)} pour tout $A \in \mathcal{F}$,
$$
\limsup_{n \to \infty} \frac{1}{v_n} \log \mu_n(A) \leqslant \sup_{\smash{\overline{A}}} s
$$

De manière générale, \textsf{(BS)} n'est pas vérifiée pour tous les mesurables. Si $\mathcal{P}$ désigne un ensemble de parties mesurables de $X$, on définit la version restreinte de la borne supérieure suivante :

\medskip

\textsf{(BS$_\mathcal{P}$)} pour tout $A \in \mathcal{P}$,
$$
\limsup_{n \to \infty} \frac{1}{v_n} \log \mu_n(A) \leqslant \sup_{\smash{\overline{A}}} s
$$

En particulier, si $D$ est une partie de $X$, on notera :

\medskip

\textsf{(BS$_{\flat , D}$)} pour tout $K \in \mathcal{F}$  tel que $K \cap D$ soit relativement compact,
$$
\limsup_{n \to \infty} \frac{1}{v_n} \log \mu_n(K) \leqslant \sup_{\smash{\overline{K}}} s
$$

Si $(\mu_n, v_n)_{n \geqslant 1}$ vérifie \textsf{(BS$_{\flat , X}$)}, on dit que $(\mu_n, v_n)_{n \geqslant 1}$ vérifie un \emph{principe de grandes déviations faible} (PGD faible). Notons que, si $D_1 \subset D_2$, alors \textsf{(BS$_{\flat , D_1}$)} entraîne \textsf{(BS$_{\flat , D_2}$)}. De plus, si $D$ est une partie de $X$ et $\nu$ une probabilité sur $\mathcal{F}$, on dira que $\nu$ est portée par $D$ si
$$
\forall (A, B) \in \mathcal{F}^2 \qquad A \cap D = B \cap D \Rightarrow \nu(A) = \nu(B)
$$
Définissons la tribu trace $\mathcal{F}|_D = \{ A \cap D \, ; \, A \in \mathcal{F} \}$. Si $\nu$ est portée par $D$, alors
$$
\nu|_D (A \cap D) = \nu(A)
$$
définit une mesure de probabilité sur $\mathcal{F}|_D$. En particulier, on voit ainsi que $\nu$ est la loi d'une variable aléatoire à valeurs dans $D$ : il suffit de considérer l'inclusion
$$
(D, \mathcal{F}|_D , \mu|_D) \hookrightarrow (X, \mathcal{F})
$$
Soit $D$ une partie de $X$ telle que, pour tout $n \geqslant 1$, $\mu_n$ soit portée par $D$. Si $(\mu_n, v_n)_{n \geqslant 1}$ vérifie \textsf{(BS$_{\flat , D}$)}, alors, sur l'espace $D$ muni de la tribu $\mathcal{F}|_D$ et de la topologie trace $\tau|_D = \{ U \cap D \, ; \, U \in \tau \}$, $(\mu_n|_D, v_n)_{n \geqslant 1}$ vérifie \textsf{(BS$_{\flat , D}$)} (sachant que, pour tout $x \in D$, $s|_D(x) = s(x)$).

\medskip

\textbf{Remarque :} Nous venons de voir que, si, pour tout $n \geqslant 1$, $\mu_n$ est portée par $D$, \textsf{(BS$_{\flat , D}$)} sur $X$ est un résultat plus fort que \textsf{(BS$_{\flat , X}$)} sur $X$ et que \textsf{(BS$_{\flat , D}$)} sur $D$. De manière générale, il n'y a pas de réciproque. Par exemple :

\medskip

$\bullet$ sur $X = \mathbb{R}$, muni de la topologie standard $\tau$ et de la tribu $\mathcal{F} = \{ \emptyset , ]-\infty , 0], ]0, +\infty[, \mathbb{R} \}$, soient, pour tout $n \geqslant 1$, $\mu_n = \delta_{(-1)^n}$ et $v_n = n$. Alors, $\underline{s} = -\infty$. Pourtant, comme tous les ensembles mesurables non vides sont non bornés, donc non relativement compacts, le PGD faible est trivialement vrai, autrement dit on a \textsf{(BS$_{\flat , X}$)} sur $X$. En revanche, si $D = \{ -1, 1 \}$, on n'a ni \textsf{(BS$_{\flat , D}$)} sur $X$, ni \textsf{(BS$_{\flat , D}$)} sur $D$.

\medskip

$\bullet$ sur $X = \mathbb{R}$, muni de la topologie standard $\tau$ et de la tribu $\mathcal{F} = \sigma(\{ ]-\infty , -2], ]-2, 0], ]0, 2], ]2, +\infty[ \}$, soient, pour tout $n \geqslant 1$, $\mu_n = \delta_{(-1)^n}$, $v_n = n$ et $D = ]-1, 1[$. Alors, $\underline{s} = -\infty$. Cette fois, \textsf{(BS$_{\flat , D}$)} sur $D$ est trivialement vérifiée : les adhérences, dans $D$, des ensembles mesurables, à savoir $]-1, 0]$, $]0, 1[$ et $]-1, 1[$, ne sont pas compactes. En revanche, on n'a ni \textsf{(BS$_{\flat , D}$)} sur $X$, ni \textsf{(BS$_{\flat , X}$)} sur $X$.

\medskip

Les deux exemples peuvent paraître pathologiques. Il n'en est rien. Autre exemple pour le premier point : sur $X = \mathbb{R}^\mathbb{R}$, muni de la tribu cylindrique $\mathcal{F}$ et de la topologie faible $\tau$, soient $x \in X \setminus \{ 0 \}$ et, pour tout $n \geqslant 1$, $\mu_n = \delta_{(-1)^n x}$. Alors, $\underline{s} = -\infty$. Pourtant, comme tous les ensembles mesurables non vides sont non bornés, donc non relativement compacts, le PGD faible est trivialement vrai, autrement dit on a \textsf{(BS$_{\flat , X}$)} sur $X$. En revanche, si $D = \{ -x, x \}$, on n'a ni \textsf{(BS$_{\flat , D}$)} sur $X$, ni \textsf{(BS$_{\flat , D}$)} sur $D$. Toutefois :

\medskip

$\bullet$ si, pour tout $n \geqslant 1$, $\mu_n$ est portée par $D$ et si les compacts de $X$ sont mesurables (en particulier si $\mathcal{F}$ est la tribu borélienne), alors
\begin{center}
\textsf{(BS$_{\flat , D}$)} sur $X$ $\iff$ \textsf{(BS$_{\flat , X}$)} sur $X$
\end{center}

$\bullet$ si, pour tout $n \geqslant 1$, $\mu_n$ est portée par $D$ et si $D$ est fermé, alors
\begin{center}
\textsf{(BS$_{\flat , D}$)} sur $X$ $\iff$ \textsf{(BS$_{\flat , D}$)} sur $D$
\end{center}

Dans le cas de Sanov, où $X = \mathcal{M}_b(E)$ est muni de la tribu cylindrique et de la topologie produit, et $D = \mathcal{M}_1^+(E)$, aucune des deux conditions n'est vérifiée.

\medskip

On peut aussi définir

\medskip

\textsf{(BS$_{c , D}$)} pour tout $C \in \mathcal{F}$ tel que $C \cap D$ soit convexe,
$$
\limsup_{n \to \infty} \frac{1}{v_n} \log \mu_n(C) \leqslant \sup_{\smash{\overline{C}}} s
$$

Dans le cas i.i.d., il existe des contre-exemples à \textsf{(BS)} déjà dans $X = \mathbb{R}^3$. En revanche, on montre que \textsf{(BS$_{c , X}$)} est vérifiée dès que $\mu_1$ est convexe-tendue (en particulier, si $X$ est un espace de Banach séparable).

\subsection{Bornes inférieures}

Si $\varphi$ est une fonction de $X$ dans $[-\infty , +\infty]$, on dit que $\varphi$ est \emph{semi-continue inférieurement} si, pour tout $t \in [-\infty , +\infty]$, l'ensemble
$$
\big\{ x \in X \, ; \, \varphi(x) > t \big\}
$$
est un ouvert de $X$. Pour toute fonction $\varphi : X \rightarrow [-\infty , +\infty]$, on note $\varphi_\bullet$ sa régularisée semi-continue inférieurement, autrement dit la plus grande fonction semi-continue inférieurement qui soit inférieure à $\varphi$.

\begin{Th}
Pour tout $x \in X$ et pour toute fonction $\varphi : X \rightarrow [-\infty , +\infty]$ mesurable, on a :
$$
\sup_{x \in X} \big( \varphi_\bullet(x) \,\textup{\d{+}}\, \underline{s}(x) \big) \leqslant \underline{p}(\varphi) \qquad \textrm{et} \qquad \sup_{x \in X} \big( \varphi_\bullet(x) \,\textup{\d{+}}\, \overline{s}(x) \big) \leqslant \overline{p}(\varphi)
$$
\end{Th}

\textbf{Démonstration :} Démontrons la première inégalité, la deuxième se montrant de façon similaire. Soient $\varphi \in \mathcal{M}$, $x \in X$, $\delta > 0$ et $M > 0$. L'ensemble
$$
V = \big\{ y \in X ; \varphi(y) \geqslant \min \big( \varphi_\bullet(x) - \delta, M \big) \big\}
$$
est mesurable et contient $\smash{\big\{ y \in X ; \varphi_\bullet(y) \geqslant \min \big( \varphi_\bullet(x) - \delta, M \big) \big\}}$, donc est un voisinage de $x$. On a alors, pour tout $n \geqslant 1$ :
\begin{align*}
\frac{1}{v_n} \log \int e^{v_n \varphi} d\mu_n & \geqslant \frac{1}{v_n} \log \int_V e^{v_n \varphi} d\mu_n\\
 & \geqslant \min \big( \varphi_\bullet(x) - \delta, M \big) + \frac{1}{v_n} \log \mu_n(V)
\end{align*}
Prenant la limite inférieure en $n$, il vient
$$
\liminf_{n \to \infty}\frac{1}{v_n} \int e^{v_n \varphi} d\mu_n \geqslant \min \big( \varphi_\bullet(x) - \delta, M \big) + \underline{s}(x)
$$
On conclut en faisant tendre $\delta$ vers $0$ et $M$ vers $+\infty$, puis en prenant le supremum en $x \in X$.\qed

\medskip

Ce résultat est ce qu'on a appelé lemme de Varadhan ouvert dans le cas i.i.d. : il est vrai en toute généralité. On peut aussi le montrer avec l'inégalité de Tchebychev dans le cas $\varphi = \lambda$ forme linéaire continue mesurable. Remarquons que la définition des fonctions $\underline{s}$ et $\overline{s}$ peut se récrire
$$
\underline{s}(x) = \inf_{A \in \mathcal{F}} \Big( \underline{p}(\delta_A) \,\dot{+}\, \big( -(\delta_A)_\bullet (x) \big) \Big) \qquad \textrm{et} \qquad \overline{s}(x) = \inf_{A \in \mathcal{F}} \Big( \overline{p}(\delta_A) \,\dot{+}\, \big( -(\delta_A)_\bullet (x) \big) \Big)
$$
où $(-\infty) \,\dot{+}\, (+\infty) = +\infty$ et, pour toute partie $A$ de $X$, on note
$$
\delta_A := -\infty \cdot 1_{X \setminus A}
$$

\begin{Cor}\label{sdef}
Pour tout $x \in X$,
$$
\underline{s}(x) = \inf_{\varphi \in \mathcal{M}} \Big( \underline{p}(\varphi) \dot{+} \big( -\varphi_\bullet(x) \big) \Big) \qquad \textrm{et} \qquad \overline{s}(x) = \inf_{\varphi \in \mathcal{M}} \Big( \overline{p}(\varphi) \dot{+} \big( -\varphi_\bullet(x) \big) \Big)
$$
\end{Cor}

Ecrites sous cette forme, les bornes inférieures font penser à la fameuse égalité
$$
s(x) = \inf_{\lambda \in X^*} \big( p(\lambda) - \lambda(x) \big) = -p^*(x)
$$
qu'on a quand $X$ est un espace vectoriel et sous certaines conditions. On notera toutefois qu'on n'a pas de lien intéressant entre entropie $s = \underline{s}$ et pression $p = \overline{p}$, pour le moment. Cela pourra se faire si $\underline{s} = \overline{s}$ ou $\underline{p} = \overline{p}$. Remarquons que, dans le cas indépendant, et même dans le cas asymptotiquement découplé (cf. \cite{CeP11}), $\underline{s} = \overline{s}$ et $\underline{p}(\lambda) = \overline{p}(\lambda)$ pour $\lambda$ forme linéraire continue et mesurable. D'autre part, la concavité de $s$ est une condition nécessaire pour obtenir l'égalité $s = -p^*$. Si $s$ n'est pas concave, le dernier corollaire suggère d'étendre l'ensemble $X^*$ à d'autres fonctions, par exemple des formes quadratiques (ce qui donne l'ensemble gaussien ; \emph{cf.} les travaux de Costeniuc, Ellis, Haven, Touchette et Turkington, par exemple \cite{EHT00} et \cite{CETT05}).

\subsection{Borne supérieure faible}

Si $\varphi$ est une fonction de $X$ dans $[-\infty , +\infty]$, on dit que $\varphi$ est \emph{semi-continue supérieurement} si, pour tout $t \in [-\infty , +\infty]$, l'ensemble
$$
\big\{ x \in X \, ; \, \varphi(x) < t \big\}
$$
est un ouvert de $X$. Pour toute fonction $\varphi : X \rightarrow [-\infty , +\infty]$, on note $\varphi^\bullet$ sa régularisée semi-continue supérieurement, autrement dit la plus petite fonction semi-continue supérieurement qui soit supérieure à $\varphi$.

\begin{Th}\label{bsf}
Soit $D$ une partie de $X$ telle que, pour tout $n \geqslant 1$, $\mu_n$ soit portée par $D$. Alors, pour tout $x \in X$ et pour toute fonction $\varphi : X \rightarrow [-\infty , +\infty]$ mesurable telle que $K = \{ \varphi > -\infty \} \cap D$ soit relativement compact, on a :
$$
\overline{p}(\varphi) \leqslant \sup_{x \in \overline{K}} \big( \varphi^\bullet(x) \,\dot{+}\, \overline{s}(x) \big)
$$
\end{Th}

\textbf{Démonstration :} Soient $\varphi \in \mathcal{M}$, $\delta > 0$ et $M > 0$. On suppose que $K = \{ \varphi > -\infty \} \cap D$ est relativement compact. Par définition de $\varphi^\bullet$, pour tout $x \in \smash{\overline{K}}$, l'ensemble
$$
U(x) = \big\{ y \in X ; \varphi(y) \leqslant \max \big( \varphi^\bullet(x) + \delta, M \big) \big\}
$$
est un voisinage mesurable de $x$. Par définition de $\overline{s}(x)$, il existe un voisinage $A(x) \in \mathcal{F}$ de $x$ tel que
$$
\limsup_{n \to \infty} \frac{1}{v_n} \log \mu_n \big( A(x) \big) \leqslant \max \big( \overline{s}(x) + \delta , M \big)
$$
Notons $V(x) = U(x) \cap A(x)$. Du recouvrement de $\smash{\overline{K}}$ par les $V(x)$ avec $x \in \smash{\overline{K}}$, on peut extraire un sous-recouvrement fini, noté $\{ V(x_i) ; i \in \{ 1, \ldots , r \} \}$. Pour tout $n \geqslant 1$, le fait que $\mu_n$ soit la loi d'une variable aléatoire à valeurs dans $D$, puis la définition de $V(x)$ donnent
\begin{align*}
\frac{1}{v_n} \log \int e^{v_n \varphi} d\mu_n
 &\leqslant \frac{1}{v_n} \log \sum_{i=1}^r \int_{V(x_i)} e^{v_n \varphi} d\mu_n \\
 &\leqslant \frac{1}{v_n} \log \sum_{i=1}^r e^{v_n \max \big( \varphi^\bullet(x) + \delta, M \big)} \mu_n \big( V(x_i) \big)
\end{align*}
Prenant la limite supérieure en $n$ et utilisant le lemme \ref{ls}, on obtient :
\begin{align*}
\limsup_{n \to \infty} \frac{1}{v_n} \log \int e^{v_n \varphi} d\mu_n
 & \leqslant \max_{1\leqslant i\leqslant r} \Big( \max \big( \varphi^\bullet(x_i) + \delta , M \big) + \max \big( (\overline{s}(x_i) + \delta) , M \big) \Big)\\
 & \leqslant \sup_{x \in \overline{K}} \Big( \max \big( \varphi^\bullet(x) + \delta , M \big) + \max \big( (\overline{s}(x) + \delta) , M \big) \Big)
\end{align*}
ce qui donne le résultat attendu, en faisant tendre $\delta$ vers $0$ et $M$ vers $+\infty$.\qed

\subsection{Condition suffisante de PGD faible}

Combinant les bornes inférieure et supérieure faible, on obtient une version du lemme de Varadhan :

\begin{Th}[Lemme de Varadhan]
Supposons $\underline{s} = \overline{s}$. Alors, pour toute fonction $\varphi : X \rightarrow ]-\infty , +\infty[$ mesurable et continue, en notant
$$
d\tilde{\mu}_n = e^{n \varphi} d\mu_n
$$
la suite $(\tilde{\mu}_n , v_n)_{n \geqslant 1}$ vérifie un PGD faible, l'entropie associée étant $\varphi + s$.
\end{Th}

\textbf{Remarque :} Nous n'avons défini la notion de PGD faible que pour des suites de mesures de probabilité. Etendre la définition à une suite de mesures quelconques ne pose pas de problème.

\medskip

En particulier, on obtient des conditions suffisantes pour qu'une suite vérifie un PGD faible. Pour tout $x \in X$, on note $\mathcal{V}_x$ un système fondamental de voisinages mesurables de $x$. On définit les propriétés (la notation est pour \og{}sous-additivité\fg{} et \og{}contrôle\fg{} du lemme \ref{lsa2}) :

\medskip

\textsf{(SAC)} Pour tout $x \in X$ et pour tout $A \in \mathcal{V}_x$,
$$
\liminf_{n \to \infty} \frac{1}{v_n} \log \mu_n(A) = \limsup_{n \to \infty} \frac{1}{v_n} \log \mu_n(A)
$$

\medskip

\textsf{(SAC$_\flat$)} Pour tout $x \in X$ et pour tout $B \in \mathcal{V}_x$, il existe $A \in \mathcal{V}_x$ tel que
$$
\limsup_{n \to \infty} \frac{1}{v_n} \log \mu_n(A) \leqslant \liminf_{n \to \infty} \frac{1}{v_n} \log \mu_n(B)
$$

\begin{Pro} \label{cspgd}
Soit $D$ une partie de $X$ telle que, pour tout $n \geqslant 1$, $\mu_n$ soit portée par $D$. On a les implications $(1) \Rightarrow (2) \Rightarrow (3) \Rightarrow (4) \Rightarrow (5) \Rightarrow (6)$ où :\\
$(1)$ \textup{\textsf{(SAC)}} ;\\
$(2)$ \textup{\textsf{(SAC$_\flat$)}} ;\\
$(3)$ $\underline{s} = \overline{s}$ ;\\
$(4)$ pour toute fonction $\varphi \in \mathcal{M}$ telle que $\{ \varphi > -\infty \} \cap D$ soit relativement compact,
$$
p(\varphi) \leqslant \sup_{x \in X} \big( \varphi^\bullet(x) \,\dot{+}\, s(x) \big)
$$
$(5)$ $(\mu_n, v_n)_{n \geqslant 1}$ vérifie \textup{\textsf{(BS$_{\flat , D}$)}} ;\\
$(6)$ $(\mu_n, v_n)_{n \geqslant 1}$ vérifie un PGD faible.
\end{Pro}

\textbf{Démonstration :} La première implication est immédiate. Pour la deuxième implication, il suffit, pour tout $x \in X$, de passer, dans \textsf{(SAC$_\flat$)}, à l'infimum en $A$ voisinage mesurable de $x$, puis en $B$ voisinage mesurable de $x$. La troisième implication est conséquence immédiate du théorème \ref{bsf}. Pour la quatrième implication, il suffit de prendre $\varphi = \delta_K$ pour $K \in \mathcal{F}$ tel que $K \cap D$ soit relativement compact. Enfin, la dernière implication a été vue dans la section Cadre et définition de cette partie.\qed

\medskip

C'est sur cette proposition que repose la théorie de Cramér pour les PGD faibles. Le lemme sous-additif du cas asymptotiquement découplé (cf. \cite{CeP11}) permet de montrer \textsf{(SAC$_\flat$)} avec les voisinage convexes mesurables. Dans le cas i.i.d. (cf. \cite{Pet11a}), on montre même la limite pour tous les convexes ouverts, et on a donc \textsf{(SAC)}. Profitons de cette remarque pour énoncer un principe de contraction. De manière générale, le principe de contraction ne s'applique pas directement aux PGD faibles car l'image réciproque, même par une application raisonnable, d'un relativement compact n'est pas relativement compacte. En revanche, la propriété \textsf{(SAC)} se transporte bien.

\begin{Th}[Principe de contraction linéaire]
Supposons que $X$ soit muni d'une structure d'espace vectoriel et que $(\mu_n, v_n)_{n \geqslant 1}$ vérifie \textup{\textsf{(SAC)}} avec, pour tout $x \in X$, $\mathcal{V}_x = \mathcal{C}_x$ l'ensemble des convexes ouverts mesurables contenant $x$. Soit $Y$ un espace vectoriel, muni d'une tribu $\mathcal{G}$ et d'une topologie $\sigma$ telles que, pour tout $y \in Y$, l'ensemble $\mathcal{C}_y$ des voisinages de $y$ convexes ouverts mesurables soit un système fondamental de voisinages. Soit $f : X \rightarrow Y$ mesurable, continue et linéaire. Alors $(\mu_n \circ f^{-1}, v_n)_{n \geqslant 1}$ vérifie \textup{\textsf{(SAC)}}. En particulier, la suite $(\mu_n \circ f^{-1}, v_n)_{n \geqslant 1}$ vérifie un PGD faible.
\end{Th}

\textbf{Démonstration :} Soient $y \in Y$ et $C \in \mathcal{C}_y$. Les hypothèses sur $f$ assurent que $f^{-1}(C)$ est un convexe ouvert mesurable de $X$. La propriété \textsf{(SAC)} permet alors de conclure que
$$
\liminf_{n \to \infty} \frac{1}{v_n} \log \mu_n \circ f^{-1} (C) = \limsup_{n \to \infty} \frac{1}{v_n} \log \mu_n \circ f^{-1} (C)
$$
Passant à l'infimum en $C \in \mathcal{C}_y$, on obtient le résultat.\qed

\subsection{Condition nécessaire et suffisante de PGD faible}

On étend ici, à notre cadre, une caractérisation des PGD faibles à l'aide de l'entropie, mentionnée dans \cite{LeP95}. Cela requiert simplement une hypothèse raisonnable de séparation avec des ensembles mesurables.

\begin{De}
Soient $X$ un ensemble, $\mathcal{F}$ une tribu sur $X$ et $\tau$ une topologie sur $X$. On dit que $X$ vérifie l'axiome de séparation $T_2^m$ si, pour tout $(x_1, x_2) \in X^2$ avec $x_1 \neq x_2$, il existe $(A_1, A_2) \in \mathcal{F}^2$ tel que $x_1 \in \smash{\aroi_1}$, $x_2 \in \smash{\aroi_2}$ et $A_1 \cap A_2 = \emptyset$.
\end{De}

On définit\footnote{On aurait pu aussi penser définir
$$
\overline{s}_{\flat, 1}(x) = \inf_{\substack{A \in \mathcal{F}\\ \aroi \ni x}} \sup_{\substack{K \in \mathcal{F}\\ K \subset\!\!\!\subset \aroi }} \limsup_{n \to \infty} \frac{1}{v_n} \log \mu_n(A) \quad \textrm{ou} \quad \overline{s}_{\flat, 2}(x) = \inf_{\substack{A \in \mathcal{F}\\ \aroi \ni x}} \sup_{\substack{K \in \mathcal{F}\\ \overline{K} \subset\!\!\!\subset \aroi }} \limsup_{n \to \infty} \frac{1}{v_n} \log \mu_n(A)
$$
mais on se rend compte que ce ne sont pas les bonnes fonctions.}
$$
\overline{s}_\flat(x) = \inf_{\substack{U \in \tau\\ U \ni x}} \sup_{\substack{K \in \mathcal{F}\\ \overline{K} \subset\!\!\!\subset U}} \limsup_{n \to \infty} \frac{1}{v_n} \log \mu_n(K)
$$
La notation $\overline{K} \subset\!\!\!\subset U$ signifie que $\overline{K}$ est un compact inclus dans $U$. On notera que $\overline{s}_\flat \leqslant \overline{s}$.



\begin{Pro}
Supposons que $X$ vérifie $T_2^m$. Alors $(\mu_n, v_n)_{n \geqslant 1}$ vérifie un PGD faible si et seulement si
$$
\overline{s}_\flat \leqslant \underline{s}
$$
\end{Pro}

\textbf{Remarque :} Sans supposer que $X$ vérifie $T_2^m$, si $(\mu_n, v_n)_{n \geqslant 1}$ vérifie un PGD faible, alors $\overline{s}_\flat \leqslant \underline{s}$.

\medskip

\textbf{Démonstration :} Pour l'implication directe, soient $A \in \mathcal{F}$ et $K \in \mathcal{F}$ tels que $\overline{K} \subset\!\!\!\subset \aroi$. Alors, le PGD faible donne
$$
\limsup_{n \to \infty} \frac{1}{v_n} \log \mu_n(K) \leqslant \sup_{\overline{K}} s \leqslant \sup_{\aroi} s \leqslant \liminf_{n \to \infty} \frac{1}{v_n} \log \mu_n(A)
$$
On en déduit, pour tout $x \in X$,
$$
\overline{s}_\flat(x) = \inf_{\substack{U \in \tau\\ U \ni x}} \sup_{\substack{K \in \mathcal{F}\\ \overline{K} \subset\!\!\!\subset U}} \limsup_{n \to \infty} \frac{1}{v_n} \log \mu_n(K) \leqslant \inf_{\substack{A \in \mathcal{F}\\ \aroi \ni x}} \liminf_{n \to \infty} \frac{1}{v_n} \log \mu_n(A) = \underline{s}(x)
$$
Pour l'autre implication, il suffit de voir que $\overline{s}_\flat$ vérifie \textsf{(BS$_\flat$)}. Pour cela, il suffit de vérifier que
$$
\mathfrak{s}(U) = \sup_{\substack{K \in \mathcal{F}\\ \overline{K} \subset\!\!\!\subset U}} \limsup_{n \to \infty} \frac{1}{v_n} \log \mu_n(K)
$$
vérifie le \emph{principle of the largest term} :
$$
\overline{\mathfrak{s}}_\flat(U_1 \cup U_2) = \overline{\mathfrak{s}}_\flat(U_1) \vee \overline{\mathfrak{s}}_\flat(U_2)
$$
En vertu du lemme \ref{ls}, cela revient à montrer que
$$
\sup_{\substack{K \in \mathcal{F}\\ \overline{K} \subset\!\!\!\subset U_1 \cup U_2}} \limsup_{n \to \infty} \frac{1}{v_n} \log \mu_n(K) = \sup_{\substack{K_1 \in \mathcal{F}\\ \overline{K_1} \subset\!\!\!\subset U_1}} \sup_{\substack{K_2 \in \mathcal{F}\\ \overline{K_2} \subset\!\!\!\subset U_2}} \limsup_{n \to \infty} \frac{1}{v_n} \log \mu_n(K_1 \cup K_2)
$$
Si $\overline{K} \subset\!\!\!\subset U_1 \cup U_2$, $\overline{K} \setminus U_2$ et $\overline{K} \setminus U_1$ sont deux compacts disjoints. L'hypothèse $T_2^m$ donne alors l'existence de $A_1, A_2 \in \mathcal{F}$ tels que $A_1 \cap A_2 = \emptyset$, $\overline{K} \setminus U_2 \subset \aroi_1$ et $\overline{K} \setminus U_2 \subset \aroi_2$. Posant $K_1 = K \setminus A_2$ et $K_2 = K \setminus A_1$, on a $\overline{K_1} \subset\!\!\!\subset U_1$, $\overline{K_2} \subset\!\!\!\subset U_2$ et $K = K_1 \cup K_2$.

\begin{center}
\def\JPicScale{.7}
\input{sbemol.pst}
\end{center}

On en déduit l'égalité voulue.\qed

\medskip

Au passage, cela permet d'affiner certains résultats, par exemple la borne supérieure. La démonstration suivante est un peu différente de la première.

\begin{Pro}
Si $X$ vérifie $T_2^m$ et si $\varphi$ est une fonction mesurable majorée telle que $K = \{ \varphi > -\infty \} \cap D$ soit relativement compact, alors
$$
p(\varphi) \leqslant \sup_{x \in \overline{K}} \big( \varphi(x) + \overline{s}_\flat(x) \big)
$$
\end{Pro}

\textbf{Démonstration :} Soit $\delta > 0$. On définit, pour $i \in \mathbb{Z}$,
$$
K(\delta, i) = K \cap \{ i \cdot \delta \leqslant \varphi \leqslant (i+1) \cdot \delta \}
$$
Comme $\varphi$ est mesurable, $K(\delta, i) \in \mathcal{F}$. De plus, notons que $\overline{K(\delta, i)} \subseteq \overline{K}$ est compact. Comme $\varphi$ est majorée, il existe $i_0 \in \mathbb{Z}$ tel que
$$
\forall i > i_0 \qquad K(\delta, i) = \emptyset
$$
Fixons $M \geqslant i_0$. On montre que :
\begin{align*}
\limsup_{n\to\infty}\frac{1}{v_n} \log \int e^{n \varphi} d\mu_n &\leqslant \limsup_{n\to\infty}\frac{1}{v_n} \log \left( \sum_{|i| \leqslant M} e^{n(i+1)\cdot \delta}\mu_n\big( K(\delta, i) \big) + e^{- nM\delta} \right)\\
 &= \left( \bigvee_{|i| \leqslant M} (i+1) \cdot \delta + \overline{s}\big( K(\delta, i) \big) \right) \vee (-M\delta)\\
 &\leqslant \sup_{x \in \overline{K}} \big( \varphi(x) + \delta + \overline{s}_\flat(x) \big) \vee (- M\delta)
\end{align*}
Reste à prendre l'infimum en $M \geqslant i_0$, puis à prendre l'infimum en $\delta > 0$.\qed

\medskip

Notons que la condition $\varphi$ majorée est en fait relative à la suite $(\mu_n, v_n)_{n \geqslant 1}$ : il suffit que
$$
\liminf_{a \to + \infty} \limsup_{n\to\infty}\frac{1}{v_n} \log \int_{|\varphi| \geqslant a} e^{nf} d\mu_n = - \infty
$$

En découle une autre version du lemme de Varadhan :

\begin{Th}[Lemme de Varadhan]
Supposons que $X$ vérifie $T_2^m$ et que $(\mu_n , v_n)_{n \geqslant 1}$ vérifie un PGD faible. Alors, pour toute fonction $\varphi : X \rightarrow ]-\infty , +\infty[$ mesurable, semi-continue inférieurement et majorée, en notant
$$
d\tilde{\mu}_n = e^{n \varphi} d\mu_n
$$
la suite $(\tilde{\mu}_n , v_n)_{n \geqslant 1}$ vérifie un PGD faible, l'entropie associée étant $\varphi + s$.
\end{Th}

\section{PGD faibles et limites projectives}

On introduit enfin une condition reliant tribu et topologie. Non seulement cette hypothèse permet de montrer que l'entropie est semi-continue supérieurement, mais elle se comporte bien vis à vis des limites projectives que nous introduisons ensuite.

\subsection{Espaces topologiques localement mesurables}

Soient $X$ un ensemble, $\mathcal{F}$ une tribu sur $X$ et $\tau$ une topologie sur $X$. On suppose que

\medskip

\textsf{(ETLM)} tout point $x$ de $X$ admet un système fondamental de voisinages mesurables $\mathcal{V}_x$.

\medskip

On dit alors que $(X, (\mathcal{V}_x)_{x \in X}, \mathcal{F}, \tau)$ est un \emph{espace topologique localement mesurable}.

\begin{Pro}
Sur un espace topologique localement mesurable, les fonctions $\underline{s}$ et $\overline{s}$ sont semi-continues supérieurement et ont pour expressions, pour tout $x \in X$,
$$
\underline{s}(x) = \inf_{V \in \mathcal{V}_x} \liminf_{n \to \infty} \frac{1}{v_n} \log \mu_n(V) \qquad \textrm{et} \qquad \overline{s}(x) = \inf_{V \in \mathcal{V}_x} \limsup_{n \to \infty} \frac{1}{v_n} \log \mu_n(V)
$$
\end{Pro}

\textbf{Démonstration :} Les expressions données de $\underline{s}$ et $\overline{s}$ découlent de la croissance de $\mu_n$, pour tout $n$. Montrons que $\underline{s}$ est semi-continue supérieurement. Soient $t \in \mathbb{R}$ et $x \in X$ tels que $s(x) < t$. Par définition de $s(x)$, il existe $V \in \mathcal{V}_x$ tel que
$$
\liminf_{n \geqslant 1} \frac{1}{v_n} \log \mu_n(V) < t
$$
Alors, pour tout $y \in V$, il existe $W \in \mathcal{V}_y$ tel que $W \subset V$, et ainsi
$$
s(y) \leqslant \liminf_{n \geqslant 1} \frac{1}{v_n} \log \mu_n(W) \leqslant \liminf_{n \geqslant 1} \frac{1}{v_n} \log \mu_n(V) < t
$$
Donc $\{ x \in X \, ; \, s(x) < t \}$ est ouvert. La démonstration est analogue pour $\overline{s}$.\qed

\subsection{Limites projectives d'espaces topologiques localement mesurables}

Soient $X$ un ensemble muni d'une tribu $\mathcal{F}$ et
$$
\overleftarrow{\mathcal{X}} = \big( X_i, f_i, f_{ij} \big)_{i \leqslant j}
$$
une famille telle que

\medskip

\textsf{(PROJ$_1$)} les indices $i$ et $j$ décrivent un ensemble $(J, \leqslant)$ préordonné filtrant à droite ;

\medskip

\textsf{(PROJ$_2$)} pour tout $i \in J$, $(X_i, (\mathcal{V}_{x, i})_{x \in X}, \mathcal{F}_i, \tau_i)$ est un espace topologique localement mesurable, $f_i$ une application mesurable de $X$ dans $X_i$ et, pour tout $(i, j) \in J^2$ tel que $i \leqslant j$, $f_{ij}$ est une application mesurable de $X_j$ dans $X_i$ telle que, pour tous $x \in X_j$ et $V \in \mathcal{V}_{f_{ij}(x), i}$,
$$
f_{ij}^{-1}(V) \in \mathcal{V}_{x, j}
$$
(en particulier, $f_{ij}$ est continue) ;

\medskip

\textsf{(PROJ$_3$)} pour tout $(i, j, k) \in J^3$ tel que $i \leqslant j \leqslant k$, on a $f_{ii} = id_{X_i}$,
$$
f_i = f_{ij} \circ f_j  \qquad \textrm{et} \qquad f_{ik} = f_{ij} \circ f_{jk}
$$

Pour tout $x \in X$, on définit
$$
\mathcal{V}_x = \{ f_i^{-1}(V_i) \, ; \, i \in J, \, V_i \in \mathcal{V}_{f_i(x), i} \}
$$
et $\tau$ la topologie initiale sur $X$ pour la famille $(f_i)_{i \in J}$. Alors, $(X, (\mathcal{V}_x)_{x \in X}, \mathcal{F}, \tau)$ est un espace topologique localement mesurable. On dit que $\overleftarrow{\mathcal{X}}$ est un \emph{système projectif d'espaces topologiques localement mesurables} et que $X$ est sa \emph{limite projective}\footnote{Au niveau des ensembles, cette définition de limite projective d'ensembles est un peu plus générale que celle de \cite[III.51]{BouE}, les structures d'espaces mesurables ne se comportant pas bien vis à vis des limites projectives.}.

\begin{Th}
Soient $\smash{\overleftarrow{\mathcal{X}} = (X_i, f_i, f_{ij})_{(i, j) \in J^2}}$ un système projectif d'espaces topologiques localement mesurables et $(X, (\mathcal{V}_x)_{x \in X}, \mathcal{F}, \tau)$ sa limite projective. Soient $(\mu_n)_{n \geqslant 1}$ une suite de mesures de probabilité sur $\mathcal{F}$ et $(v_n)_{n \geqslant 1}$ une suite de réels strictement positifs.\\
$\bullet$ Si, pour tout $i \in J$, $(\mu_n \circ f_i^{-1}, v_n)_{n \geqslant 1}$ vérifie \textup{\textsf{(SA$_\flat$)}}, alors $(\mu_n, v_n)_{n \geqslant 1}$ vérifie \textup{\textsf{(SA$_\flat$)}}.\\
$\bullet$ Pour tout $(i, j) \in J^2$, si $i \leqslant j$, on a :
$$
\underline{s}_i \circ f_i \geqslant \underline{s}_j \circ f_j \qquad \textrm{et} \qquad \overline{s}_i \circ f_i \geqslant \overline{s}_j \circ f_j
$$
De plus,
$$
\underline{s} = \inf_{i \in J} (\underline{s}_i \circ f_i) \qquad \textrm{et} \qquad \overline{s} = \inf_{i \in J} (\overline{s}_i \circ f_i)
$$
En particulier, si, pour tout $i \in J$, $\underline{s}_i = \overline{s}_i$, alors $\underline{s} = \overline{s}$.\\
$\bullet$ Soit $D$ une partie de $X$ telle que, pour tout $n \geqslant 1$, $\mu_n$ soit portée par $D$. Si, pour tout $i \in J$, les compacts de $X_i$ sont mesurables et $(\mu_n \circ f_i^{-1}, v_n)_{n \geqslant 1}$ vérifie un PGD faible, alors $(\mu_n, v_n)_{n \geqslant 1}$ vérifie \textup{\textsf{(BS$_{\flat , D}$)}}.
\end{Th}

\textbf{Démonstration :} Pour le premier point, soient $x \in X$ et $B \in \mathcal{F}$ un voisinage de $x$. Soit $V \in \mathcal{V}_x$ tel que $V \subset B$. Il existe $i \in J$ et $V_i \in \mathcal{V}_{f_i(x), i}$ tels que $V = f_i^{-1}(V_i)$. L'hypothèse \textsf{(SA$_\flat$)} appliquée à $(\mu_n \circ f_i^{-1}, v_n)_{n \geqslant 1}$ donne l'existence d'un voisinage $A_i \in \mathcal{F}_i$ tel que
$$
\limsup_{n \to \infty} \frac{1}{v_n} \log \mu_n \circ f_i^{-1}(A_i) \leqslant \liminf_{n \to \infty} \frac{1}{v_n} \log \mu_n \circ f_i^{-1} (V_i)
$$
Enfin, il existe $W_i \in \mathcal{V}_{f_i(x), i}$ tel que $W_i \subset A_i$, de sorte que $A = f_i^{-1}(W_i) \in \mathcal{V}_x$ vérifie
$$
\limsup_{n \to \infty} \frac{1}{v_n} \log \mu_n(A) \leqslant \liminf_{n \to \infty} \frac{1}{v_n} \log \mu_n(B)
$$
Pour le deuxième point, nous ferons la démonstration pour les entropies inférieures. Soit $(i, j) \in J^2$ avec $i \leqslant j$. On a
$$
\underline{s}_i \circ f_i = \underline{s}_i \circ f_{ij} \circ f_j
$$
Pour montrer le premier résultat, il suffit donc de montrer que $\underline{s}_i \circ f_{ij} \geqslant \underline{s}_j$. Et, en effet, pour tout $x \in X_j$,
\begin{align*}
\inf_{V_i \in \mathcal{V}_{f_{ij}(x), i}} \liminf_{n \to \infty} \frac{1}{v_n} \log \mu_n \circ f_i^{-1}(V_i) &= \inf_{V_i \in \mathcal{V}_{f_{ij}(x), i}} \liminf_{n \to \infty} \frac{1}{v_n} \log \mu_n \circ f_j^{-1} \circ f_{ij}^{-1}(V_i)\\
 &\geqslant \inf_{V_j \in \mathcal{V}_{x, j}} \liminf_{n \to \infty} \frac{1}{v_n} \log \mu_n \circ f_j^{-1} (V_j)
\end{align*}
car $f_{ij}^{-1}(V_i) \in \mathcal{V}_{x, j}$. Puis, pour tout $x \in X$, on a :
\begin{align*}
\underline{s}(x) &= \inf_{V \in \mathcal{V}_x} \liminf_{n \to \infty} \frac{1}{v_n} \log \mu_n(V)\\
 &= \inf_{i \in J} \inf_{V_i \in \mathcal{V}_{f_i(x), i}} \liminf_{n \to \infty} \frac{1}{v_n} \log \mu_n \big( f_i^{-1} (V_i) \big) = \inf_{i \in J} \underline{s}_i\big( f_i(x) \big)
\end{align*}

Pour le troisième point, soient $K \in \mathcal{F}$ tel que $K \cap D$ soit relativement compact et $i \in J$. Etant donné que, pour tout $n \geqslant 1$, $\mu_n$ est la loi d'une variable aléatoire à valeurs dans $D$, on a :
$$
\limsup_{n \to \infty} \frac{1}{v_n} \log \mu_n(K) \leqslant \limsup_{n \to \infty} \frac{1}{v_n} \log (\mu_n \circ f_i^{-1}) \big( f_i(\overline{K \cap D}) \big) \leqslant \sup_{ f_i( \overline{K \cap D} )  } s_i = \sup_{\overline{K \cap D}} s_i \circ f_i
$$
Il reste à montrer que
$$
\inf_{i \in J} \sup_{\overline{K \cap D}} s_i \circ f_i \leqslant \sup_{\overline{K \cap D}} \inf_{i \in J} s_i \circ f_i
$$
Soient $\delta > 0$ et $M < 0$. Le point précédent permet de voir que, pour tout $x \in \overline{K \cap D}$, il existe $i(x) \in J$ tel que
$$
s_{i(x)} \circ f_{i(x)} (x) \leqslant (s(x) + \delta/2) \wedge (M - \delta/2)
$$
Etant donné que $s_{i(x)} \circ f_{i(x)}$ est s.c.s., il existe $V(x) \in \mathcal{V}_x$ tel que
$$
\sup_{V(x)} s_{i(x)} \circ f_{i(x)} \leqslant (s_{i(x)} \circ f_{i(x)} (x) + \delta/2) \wedge M
$$
On en déduit que
$$
\sup_{V(x)} s_{i(x)} \circ f_{i(x)} \leqslant (s(x) + \delta) \wedge M
$$
Du recouvrement de $\overline{K \cap D}$ par les $V(x)$, pour $x \in \overline{K \cap D}$, on peut extraire un sous-recouvrement fini, noté $\{ V(x_k) \, ; \, k \in \{ 1, \ldots , r \} \}$. Comme $J$ est filtrant à droite, il existe $i \in J$ majorant $\{ i(x_1), \ldots , i(x_r) \}$. Alors, en vertu du premier point,
$$
\forall k \in \{ 1, \ldots , r \} \qquad s_i \circ f_i \leqslant s_{i(x_k)} \circ f_{i(x_k)}
$$
On en déduit
$$
\sup_{\overline{K \cap D}} s_i \circ f_i \leqslant \max_{1 \leqslant k \leqslant r} (s(x_k) + \delta) \wedge M \leqslant \sup_{\overline{K \cap D}} (s(x) + \delta) \wedge M
$$
On conclut en prenant l'infimum en $i \in J$, puis l'infimum en $\delta > 0$ et en $M < 0$.\qed

\medskip

On notera que les deux premiers points seraient encore vrais si on ne supposait que : pour tout $x \in X$, pour tout voisinage mesurable $V$ de $X$, il existe $i \in J$ et un voisinage mesurable $V_i$ de $f_i(x)$ tel que
$$
V \supset f_i^{-1}(V_i)
$$
Pour le dernier point en revanche, on se sert vraiment de la structure d'espace topologique localement mesurable. En ce qui concerne l'hypothèse que, pour tout $i \in J$, les compacts de $X_i$ sont mesurables, il suffit de supposer que : pour tout $K \in \mathcal{F}$ tel que $K \cap D$ soit relativement compact et pour tout $i \in J$, il existe un ensemble mesurable compris entre $f_i(K \cap D)$ et $f_i(\overline{K})$. La démonstration s'adapte alors très bien, en remarquant que, comme $K$ est relativement compact, $\smash{\overline{f_i(K)}} = f_i(\overline{K})$. Notons enfin que, sans structure d'espace topologique localement mesurable, seuls subsistent les résultats suivants :
$$
\underline{s}_i \circ f_i \geqslant \underline{s}_j \circ f_j \qquad \textrm{et} \qquad \overline{s}_i \circ f_i \geqslant \overline{s}_j \circ f_j
$$
et
$$
\underline{s} \leqslant \inf_{i \in J} (\underline{s}_i \circ f_i) \qquad \textrm{et} \qquad \overline{s} \leqslant \inf_{i \in J} (\overline{s}_i \circ f_i)
$$

\section{Equivalence d'ensembles et limites projectives}

Cette partie a pour but ultime d'énoncer une version linéaire du théorème de Dawson-Gärtner.

\subsection{Transformation de Fenchel-Legendre}

Soient $X$ un espace vectoriel, $\mathcal{F}$ une tribu sur $X$ et $\tau$ une topologie sur $X$. On note $X^{*m}$ l'ensemble des formes linéaires continues et mesurables sur $X$. On définit, pour tout $x \in X$,
$$
- \underline{p}^*(x) := \inf_{\lambda \in X^{*m}} \big( \underline{p}(\lambda) - \lambda(x) \big) \qquad \textrm{et} \qquad - \overline{p}^*(x) = \inf_{\lambda \in X^{*m}} \big( \overline{p}(\lambda) - \lambda(x) \big)
$$

\begin{Pro}
Les fonctions $- \underline{p}^*$ et $- \overline{p}^*$ sont concaves et semi-continues supérieurement.
\end{Pro}

\textbf{Démonstration :} Ce sont des infimums de fonctions affines continues.\qed

\medskip

Sous cette forme, les définitions de $- \underline{p}^*$ et $- \overline{p}^*$ rappellent fortement le corollaire \ref{sdef}, et on voit immédiatement que
$$
\underline{s} \leqslant - \underline{p}^* \qquad \textrm{et} \qquad \overline{s} \leqslant - \overline{p}^*
$$
La question de savoir s'il y a égalité est centrale dans la théorie de Cramér. Plus précisément, on s'intéresse à l'égalité reliant entropie et pression : l'égalité
$$
\underline{s} = -\overline{p}^*
$$
est-elle vérifiée ? Pour qu'elle le soit, plusieurs conditions doivent être réunies : les remarques précédentes montrent qu'il faut que $\underline{s} = \overline{s}$ et que $\underline{s}$ soit concave. Dans le cas asymptotiquement découplé (cf. \cite{CeP11}), ces deux premières propriétés sont conséquences du lemme sous-additif. Ensuite, il est nécessaire que l'ensemble $X^{*m}$ soit assez riche dans $\mathcal{M}$. Il s'agit là d'une hypothèse reliant la tribu et la topologie (cf. espaces vectoriels topologiques mesurables). De manière générale, ces conditions ne sont pas suffisantes. Notons que, s'il existe une partie $K$ de $X$ relativement compacte telle que, pour tout $n \geqslant 1$, $\mu_n$ soit portée par $K$, alors la partie sur la borne supérieure faible permet de conclure. Plus généralement, si la suite $(\mu_n)_{n \geqslant 1}$ est exponentiellement tendue, on a le résultat. Dans le cas a.d.i., le concept de convexe-tension est moins restrictif et plus pertinent pour aboutir au résultat.


\subsection{Espaces vectoriels topologiques mesurables}

Soient $X$ un espace vectoriel, $\mathcal{F}$ une tribu sur $X$ et $\tau$ une topologie sur $X$. On suppose que

\medskip

\textsf{(EVTM$_1$)} $(X, \tau)$ est un espace vectoriel topologique ;

\medskip

\textsf{(EVTM$_2$)} tout point $x$ de $X$ a un système fondamental de voisinages mesurables $\mathcal{V}_x$ ;

\medskip

\textsf{(EVTM$_3$)} toute forme linéaire continue est mesurable.

\medskip

On dit alors que $(X, (\mathcal{V}_x)_{x \in X}, \mathcal{F}, \tau)$ est un \emph{espace vectoriel topologique mesurable}. Si on note $X^*$ (resp. $X^{*m}$) l'ensemble des formes linéaires continues (resp. continues et mesurables) sur $X$, on a $X^{*m} = X^*$.



\medskip

En particulier, $(X, (\mathcal{V}_x)_{x \in X}, \mathcal{F}, \tau)$ est un espace topologique localement mesurable. Les deux autres hypothèses sont motivées par deux raisons. La première est qu'ainsi $\underline{p}^*$ et $\overline{p}^*$ sont des transformées de Fenchel-Legendre au sens habituel. En particulier, on a le théorème d'inversion (qui dit essentiellement que $X^*$ est suffisamment riche dans l'ensemble des fonctions convexes semi-continues inférieurement). La seconde est qu'alors les fonctions $-\underline{p}^*$ et $-\overline{p}^*$ se comportent bien vis à vis des limites projectives. Remarquons que les \emph{espaces localement convexes mesurables} introduits dans \cite{Pet11a} et \cite{CeP11} sont des espaces vectoriels topologiques mesurables. Les convexes y jouent un rôle crucial en lien avec la sous-additivité.

\medskip

\textbf{Remarque :} En revanche, les fonctions $(-\underline{s})^*$ et $(-\overline{s})^*$ définies pas la transformée de Fenchel-Legendre inverse se comportent mal vis à vis des limites projectives. On montre que
$$
(-\underline{s})^*(\lambda_i \circ f_i) = \sup_{x \in X} \inf_{j \geqslant i} \Big( \lambda_i \circ f_{ij} \big( f_j(x) \big) + \underline{s}_j \circ f_j(x) \Big)
$$
et on ne peut pas intervertir le supremum et l'infimum, en général. On mentionnera simplement qu'on saurait intervertir si $x$ décrivait un ensemble relativement compact (\emph{cf.} la démonstration du théorème relatif aux limites projectives d'espaces topologiques localement mesurables).

\subsection{Limites projectives d'espaces vectoriels topologiques mesurables}

Soient $X$ un espace vectoriel, muni d'une tribu $\mathcal{F}$, et
$$
\overleftarrow{\mathcal{X}} = \big( X_i, f_i, f_{ij} \big)_{i \leqslant j}
$$
une famille telle que

\medskip

\textsf{(PROJ$_1$)} les indices $i$ et $j$ décrivent un ensemble $(J, \leqslant)$ préordonné filtrant à droite ;

\medskip

\textsf{(PROJ$_2$)} pour tout $i \in J$, $(X_i, (\mathcal{V}_{x, i})_{x \in X}, \mathcal{F}_i, \tau_i)$ est un espace vectoriel topologique mesurable, $f_i$ une application linéaire mesurable de $X$ dans $X_i$ et, pour tout $(i, j) \in J^2$ tel que $i \leqslant j$, $f_{ij}$ est une application linéaire mesurable de $X_j$ dans $X_i$ telle que, pour tous $x \in X_j$ et $V \in \mathcal{V}_{f_{ij}(x), i}$,
$$
f_{ij}^{-1}(V) \in \mathcal{V}_{x, j}
$$
(en particulier, $f_{ij}$ est continue) ;

\medskip

\textsf{(PROJ$_3$)} pour tout $(i, j, k) \in J^3$ tel que $i \leqslant j \leqslant k$, on a $f_{ii} = id_{X_i}$,
$$
f_i = f_{ij} \circ f_j  \qquad \textrm{et} \qquad f_{ik} = f_{ij} \circ f_{jk}
$$

Pour tout $x \in X$, on définit
$$
\mathcal{V}_x = \{ f_i^{-1}(V_i) \, ; \, i \in J, \, V_i \in \mathcal{V}_{f_i(x), i} \}
$$
et $\tau$ la topologie initiale sur $X$ pour la famille $(f_i)_{i \in J}$. Alors, le lemme ci-après permet de voir que $(X, (\mathcal{V}_x)_{x \in X}, \mathcal{F}, \tau)$ est un espace vectoriel topologique mesurable. On dit que $\overleftarrow{\mathcal{X}}$ est un \emph{système projectif d'espaces vectoriels topologiques mesurables} et que $X$ est sa \emph{limite projective}.

\begin{Le}
Soient $\smash{\overleftarrow{\mathcal{X}} = (X_i, f_i, f_{ij})_{(i, j) \in J^2}}$ un système projectif d'espaces vectoriels topologiques mesurables et $(X, (\mathcal{V}_x)_{x \in X}, \mathcal{F}, \tau)$ sa limite projective. Alors
$$
X^* = \{ \lambda_i \circ f_i \, ; \, i \in J, \, \lambda_i \in X_i^* \}
$$
\end{Le}

\textbf{Démonstration :} Soit $\lambda \in X^*$. Comme $\lambda$ est continue, il existe $V \in \mathcal{V}_0$ tel que
$$
V \subset \{ x \in X \, ; \, \lambda(x) < 1 \}
$$
Il existe $i \in J$ et $V_i \in \mathcal{V}_{0, i}$ tels que $V = f_i^{-1}(V_i)$. On voit alors que
$$
\forall x \in X \qquad f_i(x) = 0 \Rightarrow \lambda(x) = 0
$$
En effet, si $\lambda(x) \neq 0$, il existe $t \in \mathbb{R}$ tel que $\lambda(tx) > 1$, donc $f_i(tx) \notin V_i$, d'où $f_{i}(x) \neq 0$. On en déduit que l'on peut passer au quotient dans
$$
\begin{array}{rcl}
x \in X & \longrightarrow & \lambda(x) \in \mathbb{R} \\
\searrow & & \nearrow \\
& f_i(x) \in X_i
\end{array}
$$
L'application quotient est une forme linéaire continue sur $X_i$.\qed

\begin{Th}
Soient $\smash{\overleftarrow{\mathcal{X}} = (X_i, f_i, f_{ij})_{(i, j) \in J^2}}$ un système projectif d'espaces vectoriels topologiques mesurables et $(X, (\mathcal{V}_x)_{x \in X}, \mathcal{F}, \tau)$ sa limite projective. Soient $(\mu_n)_{n \geqslant 1}$ une suite de mesures de probabilité sur $\mathcal{F}$ et $(v_n)_{n \geqslant 1}$ une suite de réels strictement positifs. Alors, pour tout $i \in J$, on a :
$$
\underline{p}_i(\lambda_i) = \underline{p}(\lambda_i \circ f_i) \qquad \textrm{et} \qquad \overline{p}_i(\lambda_i) = \overline{p}(\lambda_i \circ f_i)
$$
et, pour tout $(i, j) \in J^2$, si $i \leqslant j$, on a :
$$
- \underline{p}_i^* \circ f_i \geqslant - \underline{p}_j^* \circ f_j \qquad \textrm{et} \qquad - \overline{p}_i^* \circ f_i \geqslant - \overline{p}_j^* \circ f_j
$$
De plus,
$$
- \underline{p}^* = \inf_{i \in J} (- \underline{p}_i^* \circ f_i) \qquad \textrm{et} \qquad - \overline{p}^* = \inf_{i \in J} (- \overline{p}_i^* \circ f_i)
$$
En particulier, si, pour tout $i \in J$, $s_i = -p_i^*$, alors $s = -p^*$.
\end{Th}

\textbf{Démonstration :} On traitera à chaque fois le cas des pressions supérieures, l'autre cas étant analogue. Pour tout $i \in J$,
$$
\overline{p}_i(\lambda_i) = \limsup_{n \to \infty} \frac{1}{v_n} \log \int e^{v_n \lambda_i} d(\mu_n \circ f_i^{-1}) = \limsup_{n \to \infty} \frac{1}{v_n} \log \int e^{v_n \lambda_i \circ f_i} d\mu_n = \overline{p}(\lambda_i \circ f_i)
$$
Soit $(i, j) \in J^2$ avec $i \leqslant j$. Etant donné que $f_i = f_{ij} \circ f_j$, il suffit de montrer que $- \overline{p}_i^* \circ f_{ij} \geqslant - \overline{p}_j^*$. Soit $x \in X_j$. On a, en utilisant le point précédent et l'égalité $f_i = f_{ij} \circ f_j$,
\begin{align*}
- \overline{p}_i^* \circ f_{ij}(x) &= \inf_{\lambda_i \in X_i^*} \big( \overline{p}(\lambda_i \circ f_{ij} \circ f_j - \lambda_i \circ f_{ij}(x) \big)\\
 &\geqslant \inf_{\lambda_j \in X_j^*} \big( \overline{p}(\lambda_j \circ f_j) - \lambda_j(x) \big)\\
 &= - \overline{p}_j^*(x)
\end{align*}
Enfin, pour tout $x \in X$,
$$
\phantom\, - \overline{p}^*(x) = \inf_{\lambda \in X^*} \big( \overline{p}(\lambda) - \lambda(x) \big) = \inf_{i \in J} \inf_{\lambda_i \in X_i^*} \big( \overline{p}(\lambda_i \circ f_i) - \lambda_i \circ f_i(x) \big) = \inf_{i \in J} \Big( - p_i^*\big( f_i(x) \big) \Big) \,\qed
$$

Pour résumer les résultats précédents, on peut énoncer une version linéaire du théorème de Dawson-Gärtner dans le cadre des espaces vectoriels topologiques mesurables. On laisse au lecteur le soin d'adapter l'énoncé au cas où, pour tout $n \geqslant 1$, $\mu_n$ est portée par une partie $D$ de $X$.

\begin{Th}[Dawson-Gärtner linéaire] \label{dgl}
Soient $\smash{\overleftarrow{\mathcal{X}} = (X_i, f_i, f_{ij})_{(i, j) \in J^2}}$ un système projectif d'espaces vectoriels topologiques mesurables et $(X, (\mathcal{V}_x)_{x \in X}, \mathcal{F}, \tau)$ sa limite projective. Soient $(\mu_n)_{n \geqslant 1}$ une suite de mesures de probabilité sur $\mathcal{F}$ et $(v_n)_{n \geqslant 1}$ une suite de réels strictement positifs.\\
$\bullet$ Si, pour tout $i \in J$, $\underline{s}_i = \overline{s}_i$, alors $\underline{s} = \overline{s}$ et $(\mu_n, v_n)_{n \geqslant 1}$ vérifie un PGD faible.\\
$\bullet$ Si, pour tout $i \in J$, les compacts de $X_i$ sont mesurables et $(\mu_n \circ f_i^{-1}, v_n)_{n \geqslant 1}$ vérifie un PGD faible, alors $(\mu_n, v_n)_{n \geqslant 1}$ vérifie un PGD faible.\\
$\bullet$ Si, pour tout $i \in J$, $s_i = -p_i^*$, alors $s = -p^*$.
\end{Th}

\section{Annexe : les trois fondements de la théorie de Cramér}

Dans cette section, nous dégageons trois résultats généraux sur lesquels reposent la théorème de Cramér dans le cas indépendant.

\subsection{Lemme sous-additif}

Le premier résultat de ce type remonte à un article de M. Fekete \cite{Fek39}.
\begin{Le}\label{lsa1}
Soit $\big(u(n)\big)_{n\geqslant 1}$ une suite à valeurs dans $[0, +\infty]$. On suppose que

\medskip

\textup{\textsf{(SA)}} $u$ est sous-additive, \emph{i.e.}
$$
\forall m, n \geqslant 1 \qquad u(m+n) \leqslant u(m) + u(n)
$$

Alors,
$$
\liminf_{n\to\infty} \frac{u(n)}{n} = \inf_{n \geqslant 1} \frac{u(n)}{n}
$$
\end{Le}

\textbf{Démonstration :} Par sous-additivité, pour tous $d, m \geqslant 1$,
$$
\frac{u(dm)}{dm} \leqslant \frac{u(m)}{m}
$$
En faisant tendre $d$ vers $\infty$, on obtient
$$
\liminf_{n\to\infty} \frac{u(n)}{n} \leqslant \liminf_{d\to\infty} \frac{u(dm)}{dm} \leqslant \frac{u(m)}{m}
$$
d'où le résultat en passant à l'infimum en $m \geqslant 1$.\qed

\begin{Le}\label{lsa2}
Soit $\big(u(n)\big)_{n\geqslant 1}$ une suite à valeurs dans $[0, +\infty]$. On suppose que

\medskip

\textup{\textsf{(SA)}} $u$ est sous-additive ;

\medskip

\textup{\textsf{(C)}} $u$ est contrôlée, \emph{i.e.} il existe $N \geqslant 1$ tel que
$$
\forall n \geqslant N \qquad u(n) < + \infty
$$

Alors, la suite $\big( u(n)/n \big)_{n\geqslant 1}$ converge vers
$$
\inf_{n \geqslant 1} \frac{u(n)}{n}
$$
\end{Le}

\textbf{Démonstration :} Soient $n \geqslant m \geqslant N$. La division euclidienne de $n$ par $m$ s'écrit $n = mq + r$ avec $q \geqslant 1$ et $r \in \{ 0, \ldots , m-1 \}$ ; ainsi, la sous-additivité permet d'écrire
$$
u(n) = u(mq + r) = u\big( m(q-1) + m+r \big) \leqslant (q-1) u(m) + u(m+r)
$$
puis
$$
\frac{u(n)}{n} \leqslant \frac{u(m)}{m} + \frac{1}{n} \max_{0\leqslant i<m} u(m+i)
$$
D'où, comme $u$ est contrôlée, en faisant tendre $n$, puis $m$ vers $\infty$, on obtient
$$
\limsup_{n\to\infty} \frac{u(n)}{n} \leqslant \liminf_{m\to\infty} \frac{u(m)}{m}
$$
autrement dit la suite $\big( u(n)/n \big)_{n\geqslant 1}$ converge. D'après le lemme \ref{lsa1}, sa limite est
$$
\phantom \qquad\qquad\qquad\qquad\qquad\qquad\qquad\qquad\quad \inf_{n \geqslant 1} \frac{u(n)}{n} \qquad\qquad\qquad\qquad\qquad\qquad\qquad\qquad\quad \square
$$

\subsection{Interversion infimum-supremum}

Il est question ici du fameux \emph{principle of the largest term} de \cite{LPS93}.

\begin{Le}\label{ls}
Si $\big(u_1(n)\big)_{n\geqslant 1}$, \ldots , $\big(u_r(n)\big)_{n\geqslant 1}$ sont $r$ suites à valeurs dans $[0, +\infty]$, on a l'égalité
\[
\limsup_{n\to\infty} \frac{1}{n} \log \sum_{i=1}^r u_i(n) = \max_{1\leqslant i\leqslant r} \limsup_{n\to\infty} \frac{1}{n} \log u_i(n)
\]
\end{Le}

\textbf{Démonstration :} De l'encadrement (on rappelle que les suites sont positives)
\[
\max_{1\leqslant i\leqslant r} u_i(n) \leqslant \sum_{i=1}^r u_i(n) \leqslant r \max_{1\leqslant i\leqslant r} u_i(n)
\]
on déduit que
\[
\limsup_{n\to\infty} \frac{1}{n} \log \sum_{i=1}^r u_i(n) = \limsup_{n\to\infty} \frac{1}{n} \log \max_{1\leqslant i\leqslant r} u_i(n)
\]
Puis
\begin{eqnarray*}
\limsup_{n\to\infty} \frac{1}{n} \log \max_{1\leqslant i\leqslant r} u_i(n) & = & \lim_{n\to\infty} \sup_{k\geqslant n} \frac{1}{k} \log \max_{1\leqslant i\leqslant r} u_i(k)\\
 & = & \lim_{n\to\infty} \max_{1\leqslant i\leqslant r} \left[\sup_{k\geqslant n} \frac{1}{k} \log u_i(k)\right]\\
 & = & \max_{1\leqslant i\leqslant r} \lim_{n\to\infty} \left[\sup_{k\geqslant n} \frac{1}{k} \log u_i(k)\right]\\
 & = & \max_{1\leqslant i\leqslant r} \limsup_{n\to\infty} \frac{1}{n} \log u_i(n)
\end{eqnarray*}
La troisième égalité vient du fait facile à vérifier que
\[
\max : [-\infty , +\infty]^r \rightarrow [-\infty , +\infty]
\]
est continue.\qed

\subsection{Transformation de Fenchel-Legendre}

Pour une présentation plus large, on renvoie à \cite{Mor67} et \cite{Roc70}. La preuve présentée ici est reprise de \cite[proposition 12.2.]{Cer07}.
Soit $(X, \tau)$ un espace vectoriel localement convexe. Notons $X^*$ son dual topologique. Définissons, pour $f : X \rightarrow [-\infty, +\infty]$, sa \emph{transformée de Fenchel-Legendre} par
$$
\forall \lambda \in X^* \qquad f^*(\lambda) := \sup_{x \in X}  \big( \langle \lambda | x \rangle - f(\lambda) \big)
$$
La transformée de Fenchel-Legendre de $g : X^* \rightarrow [-\infty, +\infty]$ se définit de façon analogue.

Etendons l'addition à $[-\infty , +\infty ]$ via
$$
\forall a \in \mathbb{R} \qquad \pm\infty \,\textup{\d{$+$}}\, a = \pm\infty \qquad \textrm{et} \qquad -\infty \,\textup{\d{$+$}}\, (+\infty) = -\infty
$$
ainsi que la multiplication par un réel via
$$
\begin{array}{ll} \forall \alpha > 0 & \alpha \cdot (\pm\infty) = \pm\infty \\
\forall \alpha < 0 & \alpha \cdot (\pm\infty) = \mp\infty
\end{array} \qquad \textrm{et} \qquad 0 \cdot (\pm\infty) = 0
$$
On dit que $f : X \rightarrow [-\infty, +\infty]$ est \emph{convexe} si
$$
\forall (x, y) \in X^2 \quad \forall \alpha \in [0, 1] \qquad f(\alpha x + (1-\alpha)y) \leqslant \alpha f(x) \,\textup{\d{$+$}}\, (1-\alpha) f(y)
$$
On notera que la seule fonction convexe prenant la valeur $-\infty$ est la fonction constante de valeur $-\infty$. Les fonctions convexes autres que les deux constantes $\pm \infty$ sont habituellement dénommées \emph{fonctions convexes propres}. On dit que $f : X \rightarrow [-\infty, +\infty]$ est \emph{concave} si $-f$ est convexe. On dit que $q : X \rightarrow [-\infty , +\infty ]$ est \emph{affine} si $q$ est convexe et concave. Les fonctions affines à valeurs dans $[-\infty, +\infty]$ sont alors les fonctions affines habituelles à valeurs dans $]-\infty, +\infty[$ et les deux fonctions constantes $\pm \infty$.

Le seul résultat qui nous intéresse ici est le suivant :
\begin{Pro}\label{fl}
Soit $(X, \tau)$ un espace vectoriel localement convexe et $f : X \rightarrow [-\infty , +\infty ]$. Alors
$$
f^{**} = f
$$
si et seulement si $f$ est convexe et semi-continue inférieurement relativement à la topologie faible $\sigma(X, X^*)$.
\end{Pro}

\textbf{Remarque :} Plus précisément, la transformation de Fenchel-Legendre réalise une bijection des fonctions convexes $\sigma(X, X^*)$-s.c.i. sur $X$ sur les fonctions convexes $\sigma(X^*, X)$-s.c.i. sur $X^*$, de sorte que, si $f$ est convexe et $\sigma(X, X^*)$-s.c.i., sa transformée de Fenchel-Legendre $f^*$ prend le nom de fonction \emph{convexe-conjuguée} de $f$.

\textbf{Démonstration :} L'implication directe découle du fait que, pour toute fonction $g : X^* \rightarrow [-\infty , +\infty]$, la fonction $g^*$ est convexe et $\sigma(X, X^*)$-s.c.i. Montrons la réciproque. On remarque que, si $\lambda \in X^*$, $\langle \lambda | \cdot \rangle - f^*(\lambda)$ est la plus grande fonction affine (y compris les deux fonctions affines $\pm \infty$) dirigée par $\lambda$ et inférieure à $f$. Il s'agit donc de voir que $f$ est la borne supérieure de l'ensemble des fonctions affines continue plus petites que $f$ (y compris les deux fonctions affines $\pm \infty$). Si $f = \pm \infty$, le résultat est immédiat. Sinon, $f$ est une fonction convexe propre. Définissons l'épigraphe de $f$ par $epi(f) = \{ (x, t) \in X \times \mathbb{R} ; f(x) \leqslant t \}$. Le fait que $f$ soit convexe et $\sigma(X, X^*)$-s.c.i. assure que $epi(f)$ est convexe et fermé relativement à $\tau$. Soit $(x, t) \in X \times \mathbb{R} \setminus epi(f)$. Le théorème de Hahn-Banach dans l'espace localement convexe $X \times \mathbb{R}$ donne l'existence d'un hyperplan fermé séparant strictement $(x, t)$ et $epi(f)$. Si cet hyperplan n'est pas vertical, il correspond à une fonction affine plus petite que $f$. Sinon, l'hyperplan est de la forme $H \times \mathbb{R}$ où $H$ est un hyperplan affine fermé de $X$ et $f(x) = +\infty$. Soit alors $p$ une fonction affine continue sur $X$ telle que $p(x) > 0$ et $p(y) = 0$ pour tout $y \in H$. Comme $f$ est une fonction convexe propre, il existe une fonction affine continue et finie $q$ inférieure à $f$ : en effet, il existe $x \in X$ tel que $f(x) \in ]-\infty , +\infty[$ et un hyperplan fermé $H$ séparant $(x, f(x) - 1)$ de $epi(f)$ ; cet hyperplan $H$ n'est pas vertical et correspond à une fonction affine continue et finie $q$ inférieure à $f$. Ainsi, pour tout $\alpha > 0$, $q + \alpha p$ est toujours une fonction affine continue inférieure à $f$. Il suffit alors de choisir $\alpha$ tel que $q(x)+\alpha p(x) > t$.\qed

\bibliographystyle{alpha-fr}
\bibliography{../cramer}

\end{document}